# Nonarithmetic superrigid groups: Counterexamples to Platonov's conjecture

By Hyman Bass and Alexander Lubotzky*

**Abstract**

Margulis showed that "most" arithmetic groups are superrigid. Platonov conjectured, conversely, that finitely generated linear groups which are superrigid must be of "arithmetic type." We construct counterexamples to Platonov's Conjecture.

## 1. Platonov's conjecture that rigid linear groups are aritmetic

(1.1) *Representation rigid groups.* Let $\Gamma$ be a finitely generated group. By a *representation* of $\Gamma$ we mean a finite dimensional complex representation, i.e. essentially a homomorphism $\rho : \Gamma \longrightarrow \mathrm{GL}_n(\mathbb{C})$, for some $n$. We call $\Gamma$ *linear* if some such $\rho$ is faithful (i.e. injective). We call $\Gamma$ *representation rigid* if, in each dimension $n \geq 1$, $\Gamma$ admits only finitely many isomorphism classes of simple (i.e. irreducible) representations.

Platonov ([P-R, p. 437]) conjectured that if $\Gamma$ is representation rigid and linear then $\Gamma$ is of "arithmetic type" (see (1.2)(3) below). Our purpose here is to construct counterexamples to this conjecture. In fact our counterexamples are representation superrigid, in the sense that the Hochschild-Mostow completion $A(\Gamma)$ is finite dimensional (cf. [BLMM] or [L-M]).

The above terminology is justified as follows (cf. [L-M]): If $\Gamma = \langle s_1, \ldots, s_d \rangle$ is given with $d$ generators, then the map $\rho \mapsto (\rho(s_1), \ldots, \rho(s_d))$ identifies $R_n(\Gamma) = \mathrm{Hom}(\Gamma, \mathrm{GL}_n(\mathbb{C}))$ with a subset of $\mathrm{GL}_n(\mathbb{C})^d$. In fact $R_n(\Gamma)$ is easily seen to be an affine subvariety. It is invariant under the simultaneous conjugation action of $\mathrm{GL}_n(\mathbb{C})$ on $\mathrm{GL}_n(\mathbb{C})^d$. The algebraic-geometric quotient $X_n(\Gamma) = \mathrm{GL}_n(\mathbb{C})\backslash\backslash R_n(\Gamma)$ exactly parametrizes the isomorphism classes of semi-simple $n$-dimensional representations of $\Gamma$. It is sometimes called the $n$-dimensional "character variety" of $\Gamma$.

With this terminology we see that $\Gamma$ is *representation rigid if and only if all character varieties of* $\Gamma$ *are finite* (or zero-dimensional). In other words, there are no moduli for simple $\Gamma$-representations.

*Work partially supported by the US-Israel Binational Science Foundation.

(1.2) *Examples and remarks.* (1) If $\Gamma' \leq \Gamma$ is a subgroup of finite index then $\Gamma$ is representation rigid if and only if $\Gamma'$ is representation rigid (cf. [BLMM]). Call groups $\Gamma$ and $\Gamma_1$ (*abstractly*) *commensurable* if they have finite index subgroups $\Gamma' \leq \Gamma$ and $\Gamma'_1 \leq \Gamma_1$ which are isomorphic. In this case $\Gamma$ is representation rigid if and only if $\Gamma_1$ is so.

(2) Let $K$ be a finite field extension of $\mathbb{Q}$, $S$ a finite set of places containing all archimedean places, and $K(S)$ the ring of $S$-integers in $K$. Let $G$ be a linear algebraic group over $K$, and $G(K(S))$ the group of $S$-integral points in $G(K)$. Under certain general conditions, for semi-simple $G$ (see (2.1) below), the Margulis superrigidity theorem applies here, and it implies in particular that $G(K(S))$ is representation rigid.

(3) Call a group $\Gamma$ of "*arithmetic type*" if $\Gamma$ is commensurable (as in (1)) with a product

$$\prod_{i=1}^{n} G_i\Big(K_i(S_i)\Big),$$

where each factor is as in (2) above.

(4) Call $\Gamma$ of "*Golod-Shafarevich representation type*" if $\Gamma$ is residually finite (the finite index subgroups have trivial intersection) and $\rho(\Gamma)$ is finite for all representations $\rho$. Such groups are representation rigid. (More generally, $\Gamma$ is representation rigid if and only if $\rho(\Gamma)$ is so for all representations $\rho$ (cf. [BLMM] and the references therein).) On the other hand, such groups are linear if and only if they are finite.

Any torsion residually finite $\Gamma$ is of Golod-Shafarevich representation type (see Burnside's proof of the Burnside conjecture for linear groups). See [Go] for examples of finitely generated infinite residually finite $p$-groups.

(5) In [BLMM] one can find a much larger variety of nonlinear representation rigid groups.

(1.3) *The Platonov conjecture.* Platonov ([P-R, p. 437]) conjectured that:

*A rigid linear group is of arithmetic type* (in the sense of (1.2)(3) above).

An essentially equivalent version of this was posed much earlier as a question in [B, Question (10.4)].

The principal aim of this paper is to construct counterexamples to Platonov's conjecture (see (1.4) below). Ironically, the method we use is inspired by a construction that Platonov and Tavgen [P-T] invented to construct a counterexample to a conjecture of Grothendieck (see Section 4 below).

Our examples also refute Grothendieck's conjecture. They further have the following properties: They are *representation reductive* (all representations are semi-simple) and they are *representation superrigid.* ($\Gamma$ is representation superrigid if, for all representations $\rho$ of $\Gamma$, the dimension of the Zariski closure of $\rho(\Gamma)$ remains bounded.)

It is known (cf. [LM], [BLMM] and the references therein) that each of the conditions — representation reductive and representation superrigid — implies representation rigid.

Concretely, our counterexample to Platonov's Conjecture takes the following form.

(1.4) THEOREM. *Let* $\Gamma$ *be a cocompact lattice in the real rank* 1 *form* $G = F_{4(-20)}$ *of the exceptional group* $F_4$.

(a) *There is a finite index normal subgroup* $\Gamma_1$ *of* $\Gamma$, *and an infinite index subgroup* $\Lambda$ *of* $\Gamma_1 \times \Gamma_1$, *containing the diagonal, such that the inclusion* $v : \Lambda \to \Gamma_1 \times \Gamma_1$ *induces an isomorphism* $\hat{v} : \hat{\Lambda} \to \hat{\Gamma}_1 \times \hat{\Gamma}_1$ *of profinite completions.*

(b) *Any representation* $\rho : \Lambda \to \mathrm{GL}_n(\mathbb{C})$ *extends uniquely to a representation* $\Gamma_1 \times \Gamma_1 \to \mathrm{GL}_n(\mathbb{C})$.

(c) $\Lambda$ (*like* $\Gamma_1 \times \Gamma_1$) *is representation reductive and representation superrigid.*

(d) $\Lambda$ *is not isomorphic to a lattice in any product of groups* $H(k)$, *where* $H$ *is a linear algebraic group over a local* (*archimedean or non-archimedean*) *field* $k$.

The proof of Theorem 1.4 relies on the remarkable fact that $\Gamma$ simultaneously satisfies two qualitatively opposing conditions. On the one hand, $\Gamma$ is superrigid in $G$, because $G = F_{4(-20)}$ is among the real rank 1 groups for which Corlette [Cor] and Gromov-Schoen [G-S] have proved a Margulis type superrigidity theorem. This implies that the images $\rho(\Gamma)$ under representations $\rho$ are quite restricted.

On the other hand, $\Gamma$ is a hyperbolic group, in the sense of Gromov. Such groups share some important properties (small cancellation theory) with free groups, which are the farthest thing from rigid. In particular (nonelementary) hyperbolic groups admit many exotic quotient groups. A particular kind of quotient, furnished by a theorem of Ol'shanskii and Rips (see (3.2) below), permits us, using a construction inspired by Platonov-Tavgen [P-T], to construct the group $\Lambda$ of Theorem 1.4 satisfying (1.4)(a). In this we must also make use of the finiteness of $H_i(\Gamma, \mathbb{Z})$ $(i = 1, 2)$. This follows from results of Kumaresan and Vogan-Zuckerman (see [V-Z]), and it is this result that singles $F_{4(-20)}$ out from the other real rank 1 groups for which superrigidity is known.

With (1.4)(a), (1.4)(b) then follows from a remarkable theorem of Grothendieck ((4.2) below). Then (1.4)(c) and (d) follow from (1.4)(b) and the superrigidity properties of $\Gamma_1$.

It is instructive to compare our result (1.4) with earlier efforts to produce superrigid nonlattices.

Let $H$ be a connected Lie group and $\Lambda \leq H$ a subgroup. Call $\Lambda$ *Margulis superrigid in* $H$ if, given a homomorphism $\rho : \Lambda \longrightarrow G'(k)$, where $G'$ is an absolutely simple algebraic group over a local field $k$, $\rho(\Lambda)$ is Zariski dense in $G'$, and $\rho(\Lambda)$ is not contained in a compact subgroup of $G'(k)$, then $\rho$ extends uniquely to a continuous homomorphism $\rho_H : H \longrightarrow G'(k)$.

Let $G$ be a connected semi-simple real linear Lie group without compact factors, and let $\Gamma$ be an irreducible lattice in $G$. The Margulis Superrigidity Theorem ([Mar, VII, (5.6)]) says that, if real rank $(G) \geq 2$, then $\Gamma$ is Margulis superrigid in $G$.

A potential source of superrigid nonlattices is groups $\Gamma$ sandwiched between two superrigid arithmetic lattices $\Gamma_i$ $(i = 1, 2)$, $\Gamma_1 \leq \Gamma \leq \Gamma_2$. The following kinds of examples have been studied:

(a) $\Gamma_i = \mathrm{SL}_n(A_i)$, $n \geq 3$, $i = 1, 2$, where $A_1 = \mathbb{Z}$ and $A_2 =$ either $\mathbb{Z}[1/q]$ for some prime $q$, or $A_2 =$ the integers in a real quadratic extension of $\mathbb{Q}$.

(b) $\Gamma_1 = \mathrm{SL}_n(\mathbb{Z})$, $n \geq 3$, and $\Gamma_2 = \mathrm{SL}_{n+1}(\mathbb{Z})$.

In each case it can be shown that $\Gamma$ satisfies superrigidity. However it has been further shown that $\Gamma$ must be commensurable with either $\Gamma_1$ or $\Gamma_2$ ([V1], [V2], [LZ]).

Venkataramana [V2] also proved the following results, which exclude certain generalizations of (1.4) to higher rank groups.

Let $G$ be a connected semi-simple real linear Lie group without compact factors, and let $\Gamma$ be an irreducible lattice in $G$. Define $\Delta : G \longrightarrow G \times G$, $\Delta(x) = (x, x)$. Let $\Lambda \leq G \times G$ be a discrete Zariski dense subgroup containing $\Delta(\Gamma)$.

(1.5) THEOREM ([V2, Th. 1]). *If real rank* $(G) \geq 2$ *then* $\Lambda \leq G \times G$ *is Margulis superrigid.*

(1.6) THEOREM ([V2, Th. 2]). *If real rank* $(G) \geq 2$ *and if* $G/\Gamma$ *is not compact then* $\Lambda$ *is a lattice in* $G \times G$.

Suppose that $G = F_{4(-20)}$, and that $\Gamma_1$ and $\Lambda$ are as in (1.4). Then the representation superrigidity of $\Lambda$ could be deduced from (1.5) except for the fact that real rank (G) = 1. Venkataramana's proof of (1.5) uses the methods of Margulis, whereas our (completely different) proof of (1.4) uses the facts that $\Gamma_1$ is superrigid, $\hat{\Lambda} = \hat{\Gamma}_1 \times \hat{\Gamma}_1$, and also uses Grothendieck's Theorem 4.2 below.

We were greatly aided in this work by communication from E. Rips, Yu. Ol'shanskii, Armand Borel, Gregg Zuckerman, and Dick Gross, to whom we express our great appreciation.

## 2. Margulis superrigidity

(2.1) *The case of lattices in real Lie groups.* Let $G$ be a connected semi-simple algebraic $\mathbb{R}$-group such that $G(\mathbb{R})^0$ (identity component) has no compact factors. Let $\Gamma \leq G(\mathbb{R})$ be an irreducible lattice. *Margulis superrigidity* refers to the following property (cf. discussion following (1.4) above):

*Let $k$ be a local (i.e. locally compact nondiscrete) field. Let $H$ be a connected $k$-simple $k$-group. Let $\rho : \Gamma \longrightarrow H(k)$ be a homomorphism with Zariski dense image $\rho(\Gamma)$. Then either the closure $\overline{\rho(\Gamma)}$ (in the $k$-topology) is compact, or else $k$ is archimedean ($k = \mathbb{R}$ or $\mathbb{C}$) and $\rho$ extends uniquely to a $k$-epimorphism $\rho_G : G \longrightarrow H$, and a continuous homomorphism $G(\mathbb{R}) \longrightarrow H(k)$.*

Margulis ([Mar, VII, (5.9)]) proves this (and much more) when

$$\mathbb{R} - \text{rank } (G) \geq 2.$$

Margulis further shows that his superrigidity implies that $\Gamma$ is "arithmetic," in the following sense:

*There are a connected, simply connected, and semi-simple $\mathbb{Q}$-group $M$ in which $M(\mathbb{Z})$ is Zariski dense, and a surjective homomorphism $\sigma : M(\mathbb{R})^\circ \longrightarrow G(\mathbb{R})^\circ$ such that*

$$\ker(\sigma) \textit{ is compact},$$

*and*

$$\sigma\Big(M(\mathbb{Z}) \cap M(\mathbb{R})^\circ\Big) \textit{ is commensurable with } \Gamma.$$

Since $M(\mathbb{Z}) \cap \ker(\sigma)$ is finite (being discrete and compact), it follows that there is a finite index subgroup $\Gamma_M \leq M(\mathbb{Z})$, contained in $M(\mathbb{R})^\circ$, which $\sigma$ maps isomorphically to a finite index subgroup $\Gamma_G = \sigma(\Gamma_M) \leq \Gamma$. Thus $\Gamma$ is abstractly commensurable with the arithmetic group $\Gamma_M \leq M(\mathbb{Z})$. Further, it follows from Margulis ([Mar, VII, (6.6)]) that any homomorphism $\rho : \Gamma_M \longrightarrow \mathrm{GL}_n(\mathbb{C})$ extends, on a subgroup of finite index, to a unique algebraic homomorphism $\rho_M : M(\mathbb{C}) \longrightarrow \mathrm{GL}_n(\mathbb{C})$. Therefore the identity component of the Zariski closure of $\rho(\Gamma_M)$ is $\rho_M(M(\mathbb{C}))$. Consequently, $\Gamma_M$, and so also $\Gamma$, is representation reductive ($M$ is semi-simple) and representation superrigid, in the sense of (1.3) above. In particular, the identity component of the Hochschild-Mostow completion $A(\Gamma)$ is $M(\mathbb{C})$.

Another consequence of Margulis superrigidity is the following property:

(FAb) *If $\Gamma_1 \leq \Gamma$ is a finite index subgroup then $\Gamma_1^{ab} (= \Gamma_1/(\Gamma_1, \Gamma_1))$ is finite.*

(2.2) *When $\mathbb{R}$-rank $(G) = 1$.* Keep the notation of (2.1), but assume now that,

$$\mathbb{R} - \text{rank}(G) = 1.$$

Then it is well-known that, in general, Margulis superrigidity fails completely. For example the lattice $\Gamma = \mathrm{SL}_2(\mathbb{Z})$ in $G = \mathrm{SL}_2(\mathbb{R})$ is virtually a free group, hence the farthest thing from rigid. On the other hand, Margulis superrigidity has been established for the following rank 1 groups.

$$\begin{aligned} G &= \mathrm{Sp}(n,1) \qquad (n \geq 2) \\ &= \text{the group of isometries of quaternionic hyperbolic space,} \end{aligned}$$

and

$$\begin{aligned} G &= F_{4(-20)} \text{ (the real rank-1 form of } F_4) \\ &= \text{the group of isometries of the hyperbolic Cayley plane.} \end{aligned}$$

Margulis superrigidity was established in the above cases by K. Corlette [Cor], who treated the case when the local field $k$ is archimedean, and by M. Gromov and R. Schoen, who treated the case of non-archimedean $k$ ([G-S]).

## 3. Exotic quotients of hyperbolic groups; Ol'shanskii's theorem

(3.1) *Normal Subgroups*: *The contrast between* rank $\geq 2$ *and* rank 1. Keep the notation of (2.1). Then Margulis has shown ([Mar, VIII, (2.6)]):

*Assume that* $\mathbb{R}$-rank $(G) \geq 2$. *If* $N \triangleleft \Gamma$ ($N$ *a normal subgroup of* $\Gamma$) *then either* $N$ *or* $\Gamma/N$ *is finite.*

Further, one has, in "most" (and conjecturally all) of these cases, a qualitative form of the congruence subgroup theorem, which implies that the finite groups $\Gamma/N$ occurring above are a very restricted family.

Now suppose that $G = \mathrm{Sp}(n,1)(n \geq 2)$ or $F_{4(-20)}$ and let $X$ denote the corresponding hyperbolic space of which $G(\mathbb{R})$ is the group of isometries. This is a space of constant negative curvature. Let $\Gamma \leq G(\mathbb{R})$ be a uniform (i.e. cocompact) lattice. Then $\Gamma$ acts properly discontinuously on $X$ with compact quotient $\Gamma\backslash X$. It follows (cf. [G-H, I, (3.2)]) that $\Gamma$ is "quasi-isometric" to $X$, and so $\Gamma$ is a "hyperbolic group," in the sense of Gromov.

The point we wish to emphasize here is that a hyperbolic group which is not elementary (i.e. virtually cyclic) has an "abundance" of normal subgroups. This results from an extension of "small cancellation theory" (originally for free groups) to all nonelementary hyperbolic groups. From this it follows that *a nonelementary hyperbolic group has many "exotic" quotient groups.*

The particular kind of quotient needed for our construction of a counterexample to the Platonov conjecture was not in the literature, and so we asked two well-known experts, E. Rips and Yu. Ol'shanskii. Rips (oral communication) outlined a proof, and, independently, Ol'shanskii communicated a different proof.

(3.2) THEOREM (Ol'shanskii [Ol], Rips). *Let $\Gamma$ be a nonelementary hyperbolic group. Then $\Gamma$ has a quotient $H = \Gamma/N \neq \{1\}$ which is finitely presented and such that the profinite completion $\hat{H} = \{1\}$.*

## 4. Grothendieck's theorem and question

(4.1) *Representations and profinite completion.* For a group $\Gamma$ and a commutative ring $A$, let

$\mathrm{Rep}_A(\Gamma)$ = the category of representations $\rho : \Gamma \longrightarrow \mathrm{Aut}_A(E)$, where $E$ is any finitely presented $A$-module.

Let

$$u : \Gamma_1 \longrightarrow \Gamma$$

be a group homomorphism. It induces a "restriction functor"

$$\begin{aligned} u_A^* : \mathrm{Rep}_A(\Gamma) &\longrightarrow \mathrm{Rep}_A(\Gamma_1) \\ \rho &\longmapsto \rho \circ u. \end{aligned}$$

It also induces a continuous homomorphism of profinite completions

$$\hat{u} : \hat{\Gamma}_1 \longrightarrow \hat{\Gamma}.$$

Grothendieck discovered the following remarkable close connection between profinite completions and representation theory.

(4.2) THEOREM (Grothendieck, [Gr, (1.2) and (1.3)]). *Let $u : \Gamma_1 \to \Gamma$ be a homomorphism of finitely generated groups. The following conditions are equivalent*:

(a) $\hat{u} : \hat{\Gamma}_1 \to \hat{\Gamma}$ *is an isomorphism.*

(b) $u_A^* :\ \mathrm{Rep}_A(\Gamma) \to\ \mathrm{Rep}_A(\Gamma_1)$ *is an isomorphism of categories for all commutative rings $A$.*

(b′) $u_A^*$ *is an isomorphism of categories for some commutative ring $A \neq \{0\}$.*

(4.3) COROLLARY. *If $\hat{u}$ is an isomorphism then any properties defined in representation theoretic terms* (*like representation rigid*, *representation superrigid*, *representation reductive*, ...) *are shared by $\Gamma$ and $\Gamma_1$.*

(4.4) *Grothendieck's question.* Consider a group homomorphism

(1) $u : \Gamma_1 \longrightarrow \Gamma$

such that

(2) $\hat{u} : \hat{\Gamma}_1 \longrightarrow \hat{\Gamma}$ is an isomorphism.

Grothendieck [Gr] investigated conditions under which one could conclude that $u$ itself is an isomorphism.

If $p : \Gamma \to \hat{\Gamma}$ is the natural homomorphism, then $p : \Gamma \to p(\Gamma)$ induces an isomorphism $\hat{p} : \hat{\Gamma} \to \widehat{p(\Gamma)}$. Thus, for the above question, it is natural to assume that $\Gamma$ is *residually finite*, i.e. $p$ is injective, and likewise for $\Gamma_1$. So we shall further assume that

(3) $\Gamma_1$ and $\Gamma$ are residually finite.

Conditions (2) and (3) imply that $u$ is injective, so we can think of $\Gamma_1$ as a subgroup of $\Gamma$:

$$\Gamma_1 \leq \Gamma \leq \hat{\Gamma} = \hat{\Gamma}_1.$$

Grothendieck [Gr] indicated a large class of groups $\Gamma_1$ such that (2) and (3) imply that $u$ is an isomorphism.

In [Gr, (3.1)], Grothendieck posed the following:

*Question.* Assume (2), (3), and

(4) $\Gamma_1$ and $\Gamma$ are finitely presented.

Must $u$ then be an isomorphism?

To our knowledge this question remains open. On the other hand, if in (4), one relaxes "finitely presented" to "finitely generated," then Platonov and Tavgen [P-T] have given a counterexample. It follows from [Gru, Prop. B], that $\Gamma_1$ in the Platonov-Tavgen example is not finitely presentable. Since we also make use of their construction, we recall it below (Section 6).

For later reference we record here the following observation.

(5) Let $M \leq_{\text{fi}} \Gamma$ be a finite index subgroup, and $M_1 = u^{-1}(M) \leq \Gamma_1$. Then

$$(*) \qquad [\Gamma_1 : M_1] \leq [\Gamma : M].$$

Moreover if $\hat{u} : \hat{\Gamma}_1 \to \hat{\Gamma}$ is an isomorphism then $u|M_1$ induces an isomorphism $\hat{M}_1 \to \hat{M}$. In this case (*) is an equality.

In fact, the projection $\Gamma \to \Gamma/M$ induces an injection $\Gamma_1/M_1 \to \Gamma/M$, which is surjective if $\hat{u} : \hat{\Gamma}_1 \to \hat{\Gamma}$ is surjective. Moreover, $M \leq \Gamma$ induces an inclusion $\hat{M} \leq \hat{\Gamma}$ with the same index, and similarly for $M_1 \leq \Gamma_1$. It follows easily that if $\hat{u} : \hat{\Gamma}_1 \to \hat{\Gamma}$ is an isomorphism then so also is $\hat{M}_1 \to \hat{M}$.

Using the Ol'shanskii-Rips Theorem 3.2 and a construction of Platonov-Tavgen (see Section 6 below), we shall prove the following result in Section 7.

(7.7) Theorem. *Let $L$ be a nonelementary hyperbolic group such that $H_1(L, \mathbb{Z})$ and $H_2(L, \mathbb{Z})$ are finite. Then there is a finite index normal subgroup $L_1$ of $L$, and an infinite index subgroup $Q$ of $L_1 \times L_1$, containing the diagonal, such that the inclusion $v : Q \to L_1 \times L_1$ induces an isomorphism $\hat{v} : \hat{Q} \to \hat{L}_1 \times \hat{L}_1$.*

Then, in Section 8, we shall quote results of Kumaresan and Vogan-Zuckerman ([V-Z, Table 8.2]) that allow us to take $L$ above to be any co-compact lattice in $G = F_{4(-20)}$. In view of the Corlette-Gromov-Schoen superrigidity theorem (2.2) it follows that $L_1$ and $L_1 \times L_1$ are representation superrigid, and hence so also is $Q$, by Grothendieck's Theorem 4.2. It is easily seen that $Q$ cannot be isomorphic to a lattice in any product of archimedean and non-archimedean linear algebraic groups, and so $Q$ will be the desired counterexample to Platonov's conjecture.

## 5. G. Higman's group, and variations

A well-known construction due to G. Higman ([H], see also [S, pp. 9–10]) gives an infinite group $H$ with four generators and four relations such that $H$ has no nontrivial finite quotient groups. Higman's idea inspired Baumslag ([Baum]) to construct the group-theoretic word,

$$w(a, b) = (bab^{-1})a(bab^{-1})^{-1}a^{-2},$$

with the following remarkable property: Let $a$ and $b$ be elements of a group $L$, and $M$ a finite index normal subgroup of $L$. If $w(a, b)$ belongs to $M$ then so also does $a$.

Platonov and Tavgen used the Higman group $H$ to give a counterexample to Grothendieck's conjecture. The crucial properties of $H$ they needed were that: (a) $H$ is finitely presented; (b) $\hat{H} = \{1\}$; and (c) $H_2(H, Z) = 0$. Our counterexample to the Platonov conjecture is modeled on the Platonov-Tavgen construction. Where Platonov-Tavgen use the Higman group, as a quotient of the four generator free group, we need a group $H$ with similar properties as a quotient of a hyperbolic group $L$. To this end we need the Ol'shanskii-Rips Theorem (3.2), which furnishes many such $H$, with properties (a) and (b), as quotients of any nonelementary hyperbolic group $L$. Ol'shanskii makes clever use of the Baumslag word in his construction. If we use the Schur universal central extension of $H$ to achieve condition (c) (Section 7), this already suffices to produce an abundance of counterexamples to Grothendieck's conjecture.

In Section 6 we present the Platonov-Tavgen fiber square construction. We want to apply this to a hyperbolic superrigid lattice $L$. However, Theorem 3.2 still does not provide us with condition (c) above, the vanishing of $H_2(H, \mathbb{Z})$. To achieve this we pull back to the universal central extension of $H$. But then, to return to a group closely related to our original lattice $L$, we are obliged to make use of some cohomological finiteness properties of $L$ (see Section 7). It is these latter cohomological properties that turn out to be available only for uniform lattices in $F_{4(-20)}$. We cite the cohomology calculations needed in Section 8, then assemble all of the above ingredients in Section 9 for the proof of Theorem 1.4.

## 6. The Platonov-Tavgen fiber square

(6.1) *The fiber square.* Let $L$ be a group and $p : L \to H = L/R$ a quotient group. Form the fiber product

$$\begin{array}{ccc} P & \longrightarrow & L \\ \big\downarrow & & \big\downarrow p \\ L & \stackrel{p}{\longrightarrow} & H \end{array}$$

where

$$\begin{aligned} P &= L \times_H L \\ &= \{(x,y) \in L \times L \mid p(x) = p(y)\}. \end{aligned}$$

Let

$$u : P \longrightarrow L \times L$$

be the inclusion, and consider the diagonal,

$$\Delta : L \longrightarrow L \times L, \quad \Delta(x) = (x,x).$$

Clearly

$$\begin{aligned} P &= (R,1) \cdot \Delta L = (1,R) \cdot \Delta L \\ &\cong R \rtimes L \end{aligned}$$

(the semi-direct product with conjugation action).

(6.2) Lemma. *If $L$ is finitely generated and $H$ is finitely presented then $P$ is finitely generated.*

*Proof.* The hypotheses easily imply that $R$ is finitely generated as a normal subgroup of $L$, whence $P \cong R \rtimes L$ is finitely generated.

A result of Grunewald ([Gru, Prop. B]) suggests that $P$ is rarely finitely presented (unless $H$ is finite).

The following result is abstracted from the argument of Platonov-Tavgen [P-T]. It provides a source of counterexamples to Grothendieck's question (4.4).

(6.3) Theorem. *Assume that*:

(a) *$L$ is finitely generated*;

(b) $\hat{H} = \{1\}$; *and*

(c) $H_2(H, \mathbb{Z}) = 0$.

*Then $\hat{u} : \hat{P} \longrightarrow \hat{L} \times \hat{L}$ is an isomorphism.*

This was proved in [P-T] when $L$ is a free group and $H$ the Higman quotient (5.1). Then condition (c) follows from the fact that $H^{ab} = \{1\}$ and the fact that $H$ has a "balanced presentation" (number of generators = number of relations).

(6.4) *Some notation.* For groups $U, V$, we write $U \leq V$ (resp., $U \triangleleft V$) to indicate that $U$ is a subgroup (resp., a normal subgroup) of $V$. Putting a subscript "fi" further denotes that $U$ has finite index in $V$. For example, $U \triangleleft_{\text{fi}} V$ signifies that $U$ is a finite index normal subgroup of $V$.

When $U$ and $V$ are subgroups of a common group we write

$$\begin{aligned} Z_V(U) &= \{v \in V | vu = uv \text{ for all } u \in U\}, \text{ and} \\ Z(U) &= Z_U(U). \end{aligned}$$

We denote the integral homology groups of $U$ by

$$H_i(U) = H_i(U, \mathbb{Z}) \qquad (i \geq 0).$$

Recall that

$$H_1(U) \cong U^{ab} = U/(U,U),$$

and $H_2(U)$ is called the Schur multiplier of $U$ (see Section (7.2)).

(6.5) *Proof of* (6.3); *first steps.* Given $R \triangleleft L$ with $H = L/R$ we have the exact canonical sequences

$$1 \longrightarrow R \xrightarrow{j} L \xrightarrow{p} H \longrightarrow 1 \tag{1}$$

and

$$\hat{R} \xrightarrow{\hat{j}} \hat{L} \xrightarrow{\hat{p}} \hat{H} \longrightarrow 1. \tag{2}$$

We are interested in

$$u : P = L \times_H L \longrightarrow L \times L,$$

which we can rewrite as a homomorphism of semi-direct products,

$$u : P = R \rtimes L \longrightarrow L \rtimes L, \tag{3}$$

where $L$ acts by conjugation on both sides. In this form,

$$\hat{u} = q \rtimes Id_{\hat{L}} : \tilde{R} \rtimes \hat{L} \longrightarrow \hat{L} \rtimes \hat{L}, \tag{4}$$

where $\tilde{R}$ is the completion of $R$ in the topology induced by the profinite topology of $R \rtimes L$. In fact three natural topologies on $R$ figure here, with corresponding completions.

$$\hat{R} \twoheadrightarrow \tilde{R} \xrightarrow{q} \bar{R} \trianglelefteq \hat{L}. \tag{5}$$

The topologies are defined, respectively, by the following families of finite index normal subgroups of $R$.

$$(6)\qquad \hat{\Lambda} = \{U | U \triangleleft_{\mathrm{fi}} R\} \;\supset\; \tilde{\Lambda} = \{U \in \hat{\Lambda} | U \triangleleft L\} \;\supset\; \bar{\Lambda} = \{V \cap R | V \triangleleft_{\mathrm{fi}} L\}.$$

In each case the completion above is the inverse limit of the corresponding family of $R/U$'s.

It is clear from (4) that,

(7) $\hat{u}$ and $q : \tilde{R} \longrightarrow \hat{L}$, have the "same kernel and cokernel." In particular, $\hat{u}$ is surjective (resp., injective) if and only if $q$ is so.

Putting $C = \ker(q)$, we have an exact sequence of profinite groups,

$$(8)\qquad 1 \longrightarrow C \longrightarrow \tilde{R} \xrightarrow{q} \hat{L} \xrightarrow{\hat{p}} \hat{H} \longrightarrow 1.$$

Now Theorem 6.3 follows from (7) and Corollary 6.7 below.

(6.6) Proposition. (a) *There is a natural action of* $\hat{L}$ *on* $C$ *so that* (8) *is an exact sequence of* $\hat{L}$*-groups.*

(b) $q : \tilde{R} \longrightarrow \hat{L}$ *is surjective if and only if* $\hat{H} = \{1\}$.

*Assume now that* $\hat{H} = \{1\}$.

(c) $C$ *is central in* $\tilde{R}$ *and has trivial* $\hat{L}$*-action.*

(d) *If* $L$ *is finitely generated then there is an epimorphism*

$$\widehat{H_2(H)} \longrightarrow C,$$

*where* $H_2(H) = H_2(H, \mathbb{Z})$ *is the Schur multiplier of* $H$.

(6.7) Corollary. *If* $\hat{H} = \{1\}$ *and* $H_2(H) = 0$ *then* $q : \tilde{R} \to \hat{L}$ *is an isomorphism.*

We first give a short direct proof of Corollary 6.7, which suffices for the applications in this paper. Since there is potential interest in the fiber square construction even when we do not know that $H_2(H) = 0$, we offer also the more detailed analysis provided by Proposition 6.6.

*Proof of* (6.7). First note that (6.6)(b) follows from (6.5)(8), and so our assumption that $\hat{H} = \{1\}$ yields an exact sequence $1 \to C \to \tilde{R} \xrightarrow{q} \hat{L} \to 1$. To show that $C = \{1\}$ we must show, with the notation of (6.5)(6), that $\tilde{\Lambda} = \bar{\Lambda}$. In other words, given $U \triangleleft_{\mathrm{fi}} R, U \triangleleft L$, we must find $V \leq_{\mathrm{fi}} L$ such that $V \cap R \leq U$.

The finite group $R/U \lhd L/U$ has centralizer $Z_{L/U}(R/U) = W/U \lhd_{\mathrm{fi}} L/U$, and $(W/U) \cap (R/U) = Z(R/U)$. Since $W \lhd_{\mathrm{fi}}.\ L$ and $\hat{H} = \{1\}$, the projection $L \to H = L/R$ maps $W$ onto $H$. Thus

$$1 \to R/U \to L/U \to H \to 1$$

restricts to a finite central extension

$$1 \to Z(R/U) \to W/U \to H \to 1.$$

Since $H_2(H) = 0$ (by assumption) the latter extension splits (see (7.2) below), so we have $W/U = Z(R/U) \times (V/U)$ where $V \lhd_{\mathrm{fi}} W \lhd_{\mathrm{fi}} L$ and $V \cap R = U$. This proves (6.7).

*Proof of* (6.6). Let $e : L \longrightarrow \hat{L}$ be the canonical homomorphism. (We are not assuming that $L$ is residually finite.) For $U \in \tilde{\Lambda} = \{U | U \lhd_{\mathrm{fi}} R, U \lhd L\}$ (see (6.5)(6)), put

(1) $U^{\hat{L}}$ = the closure of $e(U)$ in $\hat{L}$,

$U^L = e^{-1}(U^{\hat{L}}) =$ "$U^{\hat{L}} \cap L$"

= the closure of $U$ in the profinite topology of $L$, and

$U^R = U^L \cap R =$ the closure of $U$ in the $\bar{\Lambda}$-topology of $R$.

(2)

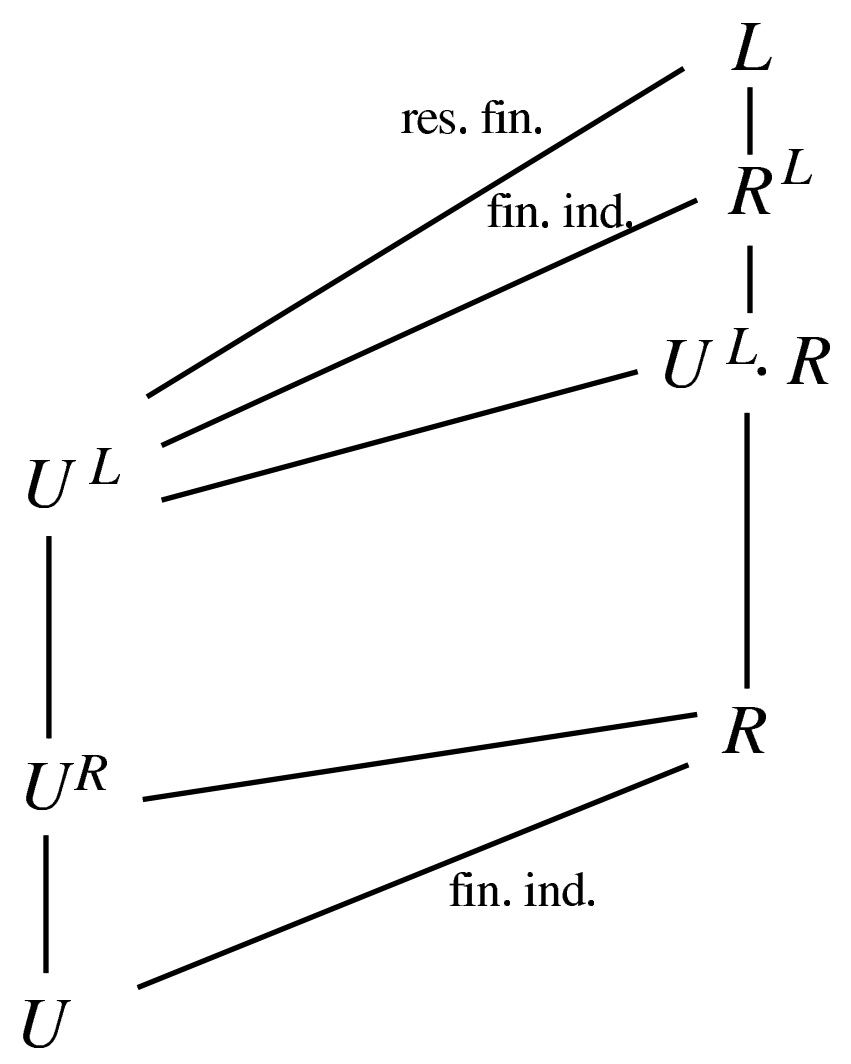

Clearly then

(3) $$U^L = \Big( \bigcap_{U \leq V \lhd_{\mathrm{fi}} L} V \Big) \lhd L,$$

and

(4) $$C(= \ker(q : \tilde{R} \longrightarrow \hat{L})) = \varprojlim_{U \in \tilde{\Lambda}} U^R/U.$$

Each $U, U^L, U^R$ above is normal in $L$, so that $C$, as an inverse limit of finite $L$-groups, hence $\hat{L}$-groups, is a profinite $\hat{L}$-group so that $q$ is equivariant for the natural action of $\hat{L}$ on $\tilde{R}$. Thus the sequence (6.5)(8) is an exact sequence of $\hat{L}$-groups, whence (6.6)(a).

For $U \in \tilde{\Lambda}$, the $L$-action on the finite group $R/U$ is continuous for the profinite $L$-topology, and $U$ acts trivially; hence $U^L$ acts trivially, i.e.

$$(U^L, R) \leq U. \tag{5}$$

Since $U^R = U^L \cap R$ it follows that

$$U^R/U \leq Z(R/U) \cap Z(U^L/U).$$

Thus,

$$(U^L \cdot R)/U \leq Z_{L/U}(U^R/U). \tag{6}$$

Next we claim that:

(7) If $\hat{H} = \{1\}$ then $U^L \cdot R = L$, $U^R/U \leq Z(L/U)$, and $L/U^L (\cong R/U^R)$ is finite.

The hypothesis implies that $R^L = L$. Since $U \triangleleft_{\text{fi}} R$ we have $U^L \triangleleft_{\text{fi}} R^L = L$. It follows that $L/U^L \cdot R$ is a finite quotient of $L/R = H$, with $\hat{H} = \{1\}$, and so $U^L \cdot R = L$. Thus $L/U^L \cong R/U^R$, and this group is finite. Finally, from (6), $U^R/U$ centralizes $U^L \cdot R/U = L/U$; i.e., $U^R/U \leq Z(L/U)$.

(8) Assume now that $\hat{H} = \{1\}$.

From (6) and (7) it follows that the conjugation induced actions of $R, \tilde{R}, L$ and $\hat{L}$ on $U^R/U$ are all trivial. Taking the inverse limit over all $U \in \tilde{\Lambda}$ we conclude from (4) that $\tilde{R}$ and $\hat{L}$ act trivially on $C$, whence (6.6)(c).

It further follows from (7) that, in the exact sequence (6.5)(1), we have $p(U^L) = H$. Thus we have from (6.5)(1) an induced central extension mod $U$,

$$1 \longrightarrow U^R/U \longrightarrow U^L/U \xrightarrow{p_1} H \longrightarrow 1. \tag{9}$$

The spectral sequence of integral homology of (9) gives an exact sequence of low order terms,

$$H_2(H) \longrightarrow U^R/U \longrightarrow H_1(U^L/U) \longrightarrow H_1(H) \longrightarrow 0. \tag{10}$$

Assume further that $L$ is finitely generated. Then so also is $H = L/R$. Since $\hat{H} = \{1\}$ it further follows that $H_1(H) = 0$.

Consider $(U^L/U, U^L/U) = (U^L, U^L) \cdot U/U$. Since $H_1(H) = 0$ we have $p_1((U^L/U, U^L/U)) = H$. The exact sequence

$$1 \longrightarrow U^L/(U^L, U^L) \cdot U \longrightarrow L/(U^L, U^L) \cdot U \longrightarrow L/U^L \longrightarrow 1$$

shows that the finitely generated group $L/(U^L, U^L) \cdot U$ is abelian-by-finite (in view of (7)), and hence residually finite. Thus $(U^L, U^L) \cdot U$ is closed in the profinite topology of $L$, whence $(U^L, U^L) \cdot U = U^L$. Thus $H_1(U^L/U) = 0$, and it follows from (10) that

(11) $$H_2(H) \longrightarrow U^R/U \text{ is surjective.}$$

Passing to the inverse limit over $U \in \tilde{\Lambda}$, we obtain a surjection

$$\widehat{H_2(H)} \longrightarrow C = \varprojlim_{U \in \tilde{\Lambda}} U^R/U,$$

whence (6.6)(d).

## 7. Schur's universal central extension

(7.1) *Application of Theorem* 6.3. Let $L$ be a nonelementary hyperbolic group. Such groups include, for example, irreducible uniform lattices in rank 1 real Lie groups. Moreover, the Corlette-Gromov-Schoen Superrigidity Theorem (see (2.2)) assures us that certain of these are representation reductive and representation superrigid.

For any nonelementary hyperbolic $L$ as above, the Ol'shanskii-Rips Theorem 3.2 furnishes us with a finitely presented quotient

(1) $$p : L \longrightarrow H = L/R \neq \{1\}, \text{ with } \hat{H} = \{1\}.$$

As in (6.1), consider the inclusion of the fiber square,

(2) $$u : P = L \times_H L \longrightarrow L \times L.$$

Since hyperbolic groups are finitely presented (see [G-H, I, (3.6)]) it follows from (6.2) that

(3) $$P \text{ is finitely generated.}$$

We have all the hypotheses of Theorem (6.3) except for

(4) $$H_2(H) = 0.$$

Given (4), we could conclude from (6.3) that $\hat{u}$ is an isomorphism. In the case that $L$ is a superrigid lattice in a Lie group $G$, this would give

$$P \quad < \quad L \times L \quad < \quad G \times G$$

with $P$ having, in view of Grothendieck's Theorem 4.2, the same representation theory as $L \times L$, hence satisfying the same superrigidity properties as $L \times L$. Yet, for these same reasons, it is easy to see that $P$ could not be a lattice in any product of real and non-archimedean linear algebraic groups. (See

Section 9 below.) Thus $P$ would furnish a strong counterexample to Platonov's conjecture 1.3.

Unfortunately this procedure assumes condition (4), and we cannot reasonably expect this to be provided by the Ol'shankii-Rips Theorem or its methods.

Instead we shall lift $p : L \to H$ to the universal central extension of $H$. We next review Schur's theory which we use for this purpose.

(7.2) *Schur's theory* (see, for example, [Mil]). Let $H$ be a group. As in (6.4), we abbreviate

$$H_i(H) = H_i(H, \mathbb{Z}) \qquad (i \geq 0).$$

We have

$$H_1(H) = H^{ab},$$

and

$$H_2(H) = \frac{R \cap (L, L)}{(L, R)} \text{ if } H = L/R, L \text{ free.}$$

(1) If $H$ is finitely presented then $H_2(H)$ is finitely generated (as a group).

Consider an exact sequence

$$1 \longrightarrow C \longrightarrow E \xrightarrow{q} H = E/C \longrightarrow 1. \tag{2}$$

(3) If $H$ and $C$ are finitely presented then so also is $E$. (See [Ha, §2, Lemma 1]).

(4) If $\hat{H} = \{1\}$, $H_1(E) = 0$, and $C$ is abelian, then $\hat{E} = \{1\}$.

In fact, if $M \triangleleft_{\text{fi}} E$ then $M \cdot C = E$ since $\hat{H} = \{1\}$, so $E/M \cong C/C \cap M$ is (like $C$) abelian, and hence trivial, since $H_1(E) = 0$.

We call (2) a *central extension* of $H$ if $C \leq Z(E)$. For any group $U$ we shall write

$$U' = (U, U).$$

(5) Suppose that (2) is central and $H_1(H) = 0$. Then $q(E') = H$ and $H_1(E') = 0$.

Clearly $q(E') = H' = H$. For $x, y \in E$, the commutator $(x, y)$ depends only on $q(x)$ and $q(y)$, since $C$ is central. We can choose $x', y' \in E'$ so that $q(x) = q(x')$ and $q(y) = q(y')$. Then $(x, y) = (x', y') \in (E', E')$, whence $E'(= (E, E)) = (E', E')$, as was to be shown.

Schur's theory tells us that:

*If* $H_1(H) = 0$ *then there is a "universal central extension,"*

$$1 \longrightarrow H_2(H) \longrightarrow \tilde{H} \xrightarrow{\pi} H \longrightarrow 1 \tag{6}$$

*which is characterized by the following equivalent properties.*

(a) (6) *is a central extension, and, given any central extension* (2), *there is a unique homomorphism* $h : \tilde{H} \to E$ *such that* $\pi = q \circ h$.

(b) (6) *is a central extension,* $H_1(\tilde{H}) = 0$, *and any central extension* $1 \to C \to E \to \tilde{H} \to 1$ *splits.*

(c) (6) *is a central extension and* $H_1(\tilde{H}) = 0 = H_2(\tilde{H})$.

(7.3) *Lifting the fiber square.* As in (6.1), consider a quotient group

$$p : L \longrightarrow H = L/R, \tag{1}$$

and the inclusion

$$u = u(p) : P = P(p) = L \times_H L \longrightarrow L \times L. \tag{2}$$

We assume that

$$\hat{H} = \{1\} \tag{3}$$

and

$$L \text{ is finitely generated.} \tag{4}$$

It follows from (3) and (4) that

$$H_1(H) = 0. \tag{5}$$

Let

$$1 \longrightarrow H_2(H) \longrightarrow \tilde{H} \xrightarrow{\pi} H \longrightarrow 1. \tag{6}$$

be the universal central extension. Form the fiber product

$$\begin{array}{ccc} \tilde{L} & \xrightarrow{\tilde{p}} & \tilde{H} \\ {\scriptstyle \pi_L}\downarrow & & \downarrow{\scriptstyle \pi} \\ L & \xrightarrow{p} & H \end{array} . \tag{7}$$

In (7), all arrows are surjective, $\ker(\pi_L) \cong \ker(\pi) = H_2(H)$, both are central, and $\ker(\tilde{p}) \cong \ker(p) = R$.

Passing to the commutator subgroups $L'$ and $\tilde{L}'$ we obtain a commutative diagram

$$\begin{array}{ccc} \tilde{L}' & \xrightarrow{\tilde{p}'} & \tilde{H} \\ {\scriptstyle \pi'_L}\downarrow & & \downarrow{\scriptstyle \pi} \\ L' & \xrightarrow{p'} & H \end{array} \tag{8}$$

in which all arrows are still surjective (since $H_1(H) = 0 = H_1(\tilde{H})$) and

(9) $\ker(\pi'_L) = \ker(\pi_L) \cap \tilde{L}'$ is isomorphic to a subgroup of $H_2(H)$.

By abelianizing the central extension

$$1 \longrightarrow \ker(\pi_L) \longrightarrow \tilde{L} \longrightarrow L \longrightarrow 1$$

one is led to the exact homology sequence of low order terms in the spectral sequence,

$$H_2(L) \longrightarrow \ker(\pi_L) \longrightarrow \tilde{L}^{ab} \longrightarrow L^{ab} \longrightarrow 0.$$

It follows from this and (9) that

(10) $\ker(\pi'_L)(= \ker(\pi_L) \cap \tilde{L}')$ is isomorphic to a quotient of $H_2(L)$.

We now assemble some conclusions.

(7.4) Proposition. *Assume that*:

(i) *$H$ is finitely presented and $\hat{H} = \{1\}$; and*

(ii) *$L$ is finitely generated and $H_1(L)$ is finite.*

*Then*:

(a) *$\tilde{H}$ is finitely presented, $\hat{\tilde{H}} = \{1\}$ and $H_2(\tilde{H}) = 0$.*

(b) *$\tilde{L}'$ is finitely generated, $L' \triangleleft_{\mathrm{fi}} L$, and there is a central extension*

$$1 \longrightarrow D \longrightarrow \tilde{L}' \xrightarrow{\pi'_L} L' \longrightarrow 1$$

*with $D$ isomorphic to a subgroup of $H_2(H)$ (and hence finitely generated) and to a quotient group of $H_2(L)$.*

(7.5) Corollary. *By* (i) *and* (ii) *of* (7.4),

$$\tilde{p}' : \tilde{L}' \longrightarrow \tilde{H}$$

*satisfies the hypotheses of Theorem* 6.3. *Hence*

$$u(\tilde{p}') : P(\tilde{p}') = \tilde{L}' \times_{\tilde{H}} \tilde{L}' \longrightarrow \tilde{L}' \times \tilde{L}'$$

*induces an isomorphism $\widehat{u(\tilde{p}')}$ of profinite completions.*

*Proof of* (7.4). The only assertion not directly covered by the discussion above is that $\tilde{L}'$ is finitely generated. Since $L$ is finitely generated (by assumption) and $H_1(L)$ is finite, $L' \triangleleft_{\mathrm{fi}} L$ is also finitely generated. Since $D$ is isomorphic to a subgroup of $H_2(H)$, which is finitely generated, it follows that $D$ is finitely generated.

To produce a counterexample to Platonov's conjecture we need conditions to insure that $\tilde{L}'$ is commensurable with $L$, or at least "close to" being so.

(7.6) *Residually finite quotients.* Let $Q$ be a group and $e : Q \longrightarrow \hat{Q}$ the canonical homomorphism to its profinite completion. We write $Q^* = e(Q)$, and call this the residually finite quotient of $Q$. Clearly $e : Q \longrightarrow Q^*$ induces an isomorphism $\hat{e} : \hat{Q} \longrightarrow \widehat{Q^*}$.

Let $u : P \longrightarrow Q$ be a group homomorphism. This induces a homomorphism $u^* : P^* \longrightarrow Q^*$, and clearly $\hat{u} : \hat{P} \longrightarrow \hat{Q}$ and $\widehat{u^*} : \widehat{P^*} \longrightarrow \widehat{Q^*}$ are isomorphic. Hence $\hat{u}$ is an isomorphism if and only if $\widehat{u^*}$ is an isomorphism.

(7.7) THEOREM. *Let $L$ be a residually finite nonelementary hyperbolic group such that $H_1(L)$ and $H_2(L)$ are finite. Then there exist $L_1 \lhd_{\rm fi} L$ and a monomorphism*

$$v : Q \longrightarrow L_1 \times L_1$$

*of infinite index such that $vQ$ contains the diagonal subgroup of $L_1 \times L_1$, and*

$$\hat{v} : \hat{Q} \longrightarrow \hat{L}_1 \times \hat{L}_1$$

*is an isomorphism.*

*Proof.* The Ol'shanskii-Rips Theorem 3.2 gives $p : L \longrightarrow H$ so that we have the hypotheses of (7.4). From (7.4) and (7.5) we then obtain the monomorphism

$$u(\tilde{p}') : P(\tilde{p}') \longrightarrow \tilde{L}' \times \tilde{L}' \tag{1}$$

of infinite index such that $\widehat{u(\tilde{p}')}$ is an isomorphism, and also a central extension

$$1 \longrightarrow D \longrightarrow \tilde{L}' \xrightarrow{\pi'_L} L' \longrightarrow 1 \tag{2}$$

with $D$ a quotient of $H_2(L)$ and hence finite, since $H_2(L)$ is assumed to be finite.

Pass to the residually finite quotient of (1).

$$u(\tilde{p}')^* : P(\tilde{p}')^* \longrightarrow \tilde{L}'^* \times \tilde{L}'^*. \tag{3}$$

Then, by (7.6), $\widehat{u(\tilde{p}')^*}$ is still an isomorphism.

Since $L$, hence also its subgroup $L' \cong \tilde{L}'/D$, is residually finite (by assumption), $\tilde{L}'^* = \tilde{L}'/D_0$ for some $D_0 \leq D$. In the residually finite group $\tilde{L}'^*$, the finite group $D/D_0$ is disjoint from some $M \leq_{\rm fi} \tilde{L}'^*$. We can even choose $M$ so that it projects isomorphically mod $D/D_0$ to a characteristic subgroup $L_1 \leq_{\rm fi} L'$. Since by our hypothesis that $H_1(L)$ is finite, $L' = (L, L) \lhd_{\rm fi} L$, it follows that $L_1 \lhd_{\rm fi} L$.

Now, in (3), $M \times M \lhd_{\rm fi} \tilde{L}'^* \times \tilde{L}'^*$ and we put $Q = u(\tilde{p}')^{*-1}(M \times M) \lhd_{\rm fi} P(\tilde{p}')^*$, and let

$$v : Q \longrightarrow M \times M$$

be the inclusion, clearly still of infinite index and containing the diagonal subgroup of $M \times M$. It follows from (4.4) (5) that $\hat{v}$ is an isomorphism. Since $M \cong L_1 \triangleleft_{\mathrm{fi}} L$, this completes the proof of (7.7).

## 8. Vanishing second Betti numbers

(8.1) Let $G$ be a connected semi-simple real Lie group with finite center, and $\Gamma \leq G$ a uniform (cocompact) lattice. The Betti numbers of $\Gamma$ are

$$b_i(\Gamma) = \dim_{\mathbb{R}} H^i(\Gamma, \mathbb{R}).$$

Since $\Gamma$ is virtually of type (FL) ([Br, VIII, 9, Ex. 4]),

$$H_i(\Gamma, \mathbb{Z}) \text{ is a finitely generated } \mathbb{Z}\text{-module},$$

and

$$H^i(\Gamma, \mathbb{R}) \cong \mathrm{Hom}_{\mathbb{R}}(H_i(\Gamma, \mathbb{Z}), \mathbb{R}).$$

Thus,

$$H_i(\Gamma, \mathbb{Z}) \text{ is finite if and only if } b_i(\Gamma) = 0.$$

(8.2) *Rank* 1 *groups.* We are interested in the case where $G$ is one of the rank 1 groups — $\mathrm{Sp}(n,1), n \geq 2$, or $F_{4(-20)}$ — for which $\Gamma$ is a hyperbolic group, and we have the Corlette-Gromov-Schoen Superrigidity Theorem 2.2. In this case $b_1(\Gamma) = 0$ because of superrigidity. We also need $b_2(\Gamma) = 0$ in order to apply Theorem 7.7.

According to Kumaresan and Vogan and Zuckerman, [V-Z], this is the case for $F_{4(-20)}$.

(8.3) Theorem ([V-Z, Table 8.2]). *If $\Gamma$ is a uniform lattice in $F_{4(-20)}$ then*

$$b_1(\Gamma) = b_2(\Gamma) = b_3(\Gamma) = 0.$$

In fact, as pointed out to us by Dick Gross, these are essentially the only examples. For, Gross indicated that it follows from result of J.-S. Li ([Li, Cor. (6.5)]), that if $\Gamma$ is a uniform lattice in $\mathrm{Sp}(n,1)$ then $b_2(\Gamma_1) \neq 0$ for some $\Gamma_1 \leq_{\mathrm{fi}} \Gamma$.

## 9. Proof of Theorem 1.4

Let $L_1$ be a cocompact lattice in $G = F_{4(-20)}$

(1) $L_1$ is a nonelementary hyperbolic group (cf. (3.1)).

(2) $H_1(L_1)$ and $H_2(L_1)$ are finite (Theorem 8.3).

Hence:

(3) There exists a finite index normal subgroup $L \triangleleft_{\rm fi} L_1$, and an infinite index subgroup $Q < L \times L$, containing the diagonal subgroup of $L \times L$ such that $\hat{Q} \longrightarrow \hat{L} \times \hat{L}$ is an isomorphism (Theorem 7.7).

Hence:

(4) Any representation $\rho : Q \longrightarrow \mathrm{GL}_n(\mathbb{C})$ extends uniquely to a representation $\rho : L \times L \longrightarrow \mathrm{GL}_n(\mathbb{C})$ (Theorem 4.2).

(5) $L \times L$, hence also $Q$ (by (4)), is representation reductive and superrigid in $G \times G$ (cf. (2.1)).

To conclude the proof of Theorem 1.4 we establish:

(6) $Q$ is not isomorphic to a lattice in any product of linear algebraic groups over archimedean and non-archimedean fields.

To prove (6), suppose that $Q$ is embedded as a lattice in $H = H_1 \times \ldots \times H_n$, where $H_i = H_i(F_i)$ is a linear algebraic group over a local field $F_i$, and so that the image of the projection $p_i : Q \longrightarrow H_i$ is Zariski dense in $H_i$. (Otherwise replace $H_i$ by the Zariski closure of $p_i Q$.) Since the representation theory of $Q$, like that of $L \times L$ (by ((4)), is semi-simple and rigid, it follows that each $H_i$ is semi-simple, and so we can even take the $H_i$ to be (almost) simple. At the cost perhaps of factoring out a finite normal subgroup of $Q$, we can further discard any $H_i$ that is compact. We further have that the (topological) closure $K_i$ of $p_i Q$ is not compact. For, say that $K_1$ was compact. Then $H/Q$ would project to the quotient $H_1/K_1 = H/(K_1 \times H_2 \times \ldots \times H_n)$; but $H_1/K_1$ does not have finite invariant volume, and this contradicts the assumption that $Q$ is a lattice in $H$.

Now it follows from (5) that $p_i : Q \longrightarrow H_i$ extends to $p_i : L \times L \longrightarrow H_i$. Let $L_1 = (L, 1)$ and $L_2 = (1, L)$, so that $L \times L = L_1 L_2$. The Zariski closures of $p_i L_1$ and $p_i L_2$ commute and generate $H_i$ ($p_i Q$ is Zariski dense). Hence, since $H_i$ is (almost) simple, one of $p_i L_j$, say $p_i L_2$, is finite (and central). Then $p_i L_1$ is Zariski dense, with noncompact closure. Now, by the Corlette-Gromov-Schoen superrigidity Theorem 2.1, there is a continuous homomorphism $q_i : G_1 = (G, 1) \longrightarrow H_i$ agreeing with $p_i$ on a finite index subgroup $L' \leq L_1$.

If $F_i$ is non-archimedean then $H_i$ is totally disconnected , so that $q_i$ must be trivial ($G$ is connected), contradicting the fact that $p_i L_1$ is not relatively compact.

Thus $F_i$ is archimedean ($\mathbb{R}$ or $\mathbb{C}$). Since $G$ is simple, $\ker(q_i)$ is finite. Moreover $q_i G_1$, like $p_i Q$, is Zariski dense in $H_i^\circ$ (identity component of $H_i$). Since $H_i^\circ$ is simple, the adjoint representation of $H_i^\circ$ on $\mathrm{Lie}(H_i^\circ)$ is irreducible.

Hence it is also an irreducible $q_iG_1$-representation. But Lie$(q_iG_1)$ in Lie $(H_i^\circ)$ is $q_iG_1$ invariant. Hence Lie $(q_iG_1) = \text{Lie}(H_i^\circ)$, so that $q_iG_1 = H_i^\circ$ .

*Conclusion.* We can identify $H_i^\circ$ with $G/Z$, $Z$ finite, so that, on a finite index subgroup of $Q$, $p_i : Q \longrightarrow H_i$ agrees with the composite $G \times G \xrightarrow{pr_1} G \longrightarrow G/Z$.

Now this happens for each $i$, except that we may sometimes use $pr_2$ in place of $pr_1$. In fact this must happen, since otherwise $Q \cap L_2$, an infinite group, would have finite image under the map $p : Q \longrightarrow H$, which has at most finite kernel.

If, say, $p_i$ factors (virtually) through $pr_i (i = 1, 2)$, then

$$p' = (p_1, p_2) : L \times L \longrightarrow H_1 \times H_2$$

virtually agrees with the projection $G \times G \longrightarrow (G/Z_1) \times (G/Z_2)$, each $Z_i$ finite. It follows that the image of $L \times L$ in $H_1 \times H_2$ is a lattice, and the image of $Q$ has infinite index in that of $L \times L$. Hence $(H_1 \times H_2)/p'Q$ cannot have finite volume, contradicting the assumption that $pQ < H$ is a lattice.

Columbia University, New York, NY
*Current address*: University of Michigan, Ann Arbor, MI
*E-mail address*: hybass@umich.edu

Institute of Mathematics, Hebrew University, Jerusalem, Israel
*E-mail address*: alexlub@math.huji.ac.il